\documentclass[a4paper]{article}

\usepackage[english]{babel} 
\usepackage[utf8x]{inputenc} 
\usepackage{amsmath}
\usepackage{amssymb}
\usepackage{amsfonts}
\usepackage{amsthm}
\usepackage{float}
\usepackage{hyperref}
\usepackage{tikz,xcolor}
\usepackage{soul}
\usepackage[normalem]{ulem}
\usepackage{ctable}     
\usepackage{multirow}	
\usepackage{mathtools}
\usepackage{multicol}
\usepackage{graphicx}
\usepackage{enumerate}

\newcommand{\Om}{\Omega}
\newcommand{\dOm}{\partial\Omega}
\newcommand{\dd}{\partial}
\newcommand{\dbar}{\overline{\partial}}

\newcommand{\R}{\mathbb{R}}

\newcommand{\F}{{\mathcal F}}
\newcommand{\N}{{\mathcal N}}
\newcommand{\cR}{{\mathcal R}}

\graphicspath{ {./images/} }

\title{Stroke classification using Virtual Hybrid Edge Detection from \textit{in silico} electrical impedance tomography data}
\author{
Juan Pablo Agnelli\thanks{FaMAF, National University of C\'ordoba and CIEM, National Scientific and Technical Research Council (CONICET), Argentina} 
\and Fernando S. Moura\thanks{Engineering, Modeling and Applied Social Sciences Center, Federal University of ABC, S\~{a}o Paulo, Brazil}
\and Siiri Rautio\thanks{Department of Mathematics and Statistics, University of Helsinki, Finland }
\and Melody Alsaker\thanks{Department of Mathematics, Gonzaga University, Spokane, USA}
\and Rashmi Murthy\thanks{Department of Mathematics, Bangalore University, India}  
\and Matti Lassas\footnotemark[3]
\and Samuli Siltanen\footnotemark[3]}

\begin{document}
	
\maketitle
\begin{abstract}
Electrical impedance tomography (EIT) is a non-invasive imaging method for recovering the internal conductivity of a physical body from electric boundary measurements. EIT combined with machine learning has shown promise for the classification of strokes. However, most previous works have used raw EIT voltage data as network inputs. We build upon a recent development which suggested the use of special noise-robust Virtual Hybrid Edge Detection (VHED) functions as network inputs, although that work used only highly simplified and mathematically ideal models. In this work we strengthen the case for the use of EIT, and VHED functions especially, for stroke classification. We design models with high detail and mathematical realism to test the use of VHED functions as inputs. Virtual patients are created using a physically detailed 2D head model which includes features known to create challenges in real-world imaging scenarios. Conductivity values are drawn from statistically realistic distributions, and phantoms are afflicted with either hemorrhagic or ischemic strokes of various shapes and sizes. Simulated noisy EIT electrode data, generated using the realistic Complete Electrode Model (CEM) as opposed to the mathematically ideal continuum model, is processed to obtain VHED functions.  We compare the use of VHED functions as inputs against the alternative paradigm of using raw EIT voltages. Our results show that (i) stroke classification can be performed with high accuracy using 2D EIT data from physically detailed and mathematically realistic models, and (ii) in the presence of noise, VHED functions outperform raw data as network inputs.
\end{abstract}

\section{Introduction}\label{sec:intro}
Stroke is a serious medical condition in which blood supply to part of the brain is blocked or significantly decreased. This prevents oxygen and other nutrients from reaching that part of the brain. Strokes are mainly classified into two groups: ischemic and hemorrhagic \cite{Donnan2008}. 
Ischemic stroke occurs when brain blood vessels become narrowed or blocked, impeding part of the brain from receiving blood. 
Hemorrhagic stroke occurs when a blood vessel ruptures or leaks, releasing blood inside the brain. 

Brain cells start to die in minutes without proper blood supply, and therefore stroke is an emergency situation that requires quick medical intervention to reduce brain damage and improve the prognosis of the patient \cite{Hacke2004, Saver2006,Saver2010}. The treatment depends on the type of stroke. Ischemic strokes are treated with tissue plasminogen activator medication to dissolve the clots, or, in some cases, thrombectomy to mechanically remove the clots from the affected blood vessel \cite{Barthels2020,Evans2017,Hacke2008}. Hemorrhagic strokes, on the other hand, require medicine to reduce blood pressure and, if possible, a medical procedure to control the bleed, such as aneurysm clipping, coiling, or AVM endovascular embolization \cite{Cordonnier2018,Hemphill2015,Qureshi2009}. The medicine administered to treat one type of stroke causes damaging effects in the other type. That is, plasminogen activator medication (for ischemic treatment) causes reduction in coagulation capability, which is a serious issue for the hemorrhagic stroke patient, while medicine to reduce blood pressure (for hemorrhagic treatment) increases obstruction problems in the ischemic stroke patient.

Because of the antagonistic effects of the medicines, early diagnosis is key, as diagnosis must be made prior to beginning treatment. Unfortunately, ischemic and hemorrhagic strokes present very similar symptoms in visual clinical inspection \cite{Ojaghihaghighi2017}. At the time of this writing, X-ray computed tomography (CT)
and magnetic resonance imaging (MRI) are the brain imaging techniques used
to make a diagnosis, and therefore diagnosis can only be performed after transporting a patient to the hospital \cite{Chalela2007,Provost2019, Rapillo2024}. If a portable diagnostic method were developed, however, the necessary diagnostic equipment could be carried in the ambulance to permit earlier treatment, which could greatly improve patient outcomes.

Electrical impedance tomography (EIT) is a relatively new tomographic imaging technique which is a good candidate for this purpose, since EIT devices are portable and do not require radiation shielding or other bulky and expensive setups. EIT generates tomographic images of the electrical conductivity distribution on a target by measuring the electrical potential on its external surface when an electrical current is applied to the target  \cite{CHENEY-1999}. 

In this work, we combine specially constructed VHED functions processed from EIT measurements with a neural network for stroke classification, and we validate the method using numerically simulated phantoms which are much more physically detailed and mathematically realistic than those used in a previous preliminary work on this topic \cite{agnelli2020classification}. Other machine learning-based approaches have been proposed for the classification of strokes using EIT data, but these approaches have all considered learning directly from measured EIT voltage data.

In \cite{macdermott2018brain} a Support Vector Machine (SVM) was applied directly to simulated EIT voltage measurements to detect the presence or absence of a hemorrhage in the brain. A 3D model was considered, and a ring of electrodes was positioned in the exterior of the model. The highly simplified head model consisted of only 2 layers: an inner brain layer and an outer layer modeling tissues external to the brain. 
Therefore, the important effects produced by the highly resistive skull and
the highly conductive cerebrospinal fluid (CSF) were not considered.
The impact of noise, lesion location and size, electrode positioning and different anatomies of the head were considered to test the methodology.
A similar approach was presented in \cite{Culpepper2024}, but considered a more realistic 3D head model consisting of 4 layers. This model consists of a first layer that aggregates the scalp and skull layers, followed by separate layers representing the CSF, brain tissue (aggregating gray and white matter), and finally, the ventricle layer.
Additionally, the study also explored the impact of four different measurement frequencies on stroke identification and differentiation using SVM.
Multi-frequency symmetry difference EIT combined with SVM was explored in \cite{McDermott2020}. This study employed a 3D head model comprising scalp, skull, CSF, and brain layers. The methodology was applied to both simulated data and a limited set of human data (35 individuals).

In \cite{shi2022RCNN} a Residual Neural Network was proposed for the classification of strokes in different regions of the head, again using minimally processed simulated EIT voltage measurements as inputs. A simplified 2D head model consisting of three layers representing the scalp, skull, and parenchyma (brain) was used. However, the highly conductive CSF
and other anatomical tissues were not considered in this model. The performance of the method was evaluated when the data was corrupted by noise, when shape deformations occurred, and when the conductivity varied in the three layers of the model. However, these conditions were considered individually, while in real-world scenarios these variations would be simultaneous. Furthermore, in that study, difference EIT data was considered as an input to the network, so the proposed method requires healthy data that may not be available in practice. 
In \cite{candiani2022NN}, fully connected neural networks (FCNN) and convolutional neural networks (CNN) were employed for detection of brain hemorrhages from simulated absolute EIT data on a 3D head model, again using raw voltages as network inputs. The model in that work was also a simplified 3-layer version of the human head consisting of scalp, skull, and brain, again neglecting the conductive cerebrospinal fluid and other anatomical tissues. The performance of the two neural networks were tested using several sets of unseen data, and FCNN showed better generalization capabilities than CNN.

A novel methodology connecting EIT and X-ray tomography was introduced in~\cite{greenleaf2018propagation} and recently analyzed in more detail in~\cite{alsaker2024VHPT}.
Using this new methodology, EIT measurements can be processed in a non-linear manner, and relevant geometric conductivity information concealed within these measurements can be extracted.
Moreover, this information hidden within EIT data possesses the same linear geometry as data produced by parallel-beam computed tomography. These ``virtual" X-ray attenuation profiles are called Virtual Hybrid Edge Detection (VHED) functions.

VHED was connected with neural network models in \cite{agnelli2020classification}  for classifying strokes as ischemic or hemorrhagic. In that introductory work, VHED functions were computed from EIT data generated by applying the unrealistic mathematically ideal continuum model of EIT to highly simplified head phantoms on circular domains. The computational evidence presented there strongly suggests that (i) when applying machine learning to severely ill-posed problems, it pays  to extract noise-robust features from the data to use as inputs for machine learning, as opposed to learning directly from measured data, and (ii) VHED profiles are useful features for classifying stroke, at least with data simulated from highly mathematically and physically unrealistic  two-dimensional models. Although the results in that initial work demonstrated the potential advantages of using VHED functions in stroke classification, it was left to future work to determine whether the method would prove robust under more realistic conditions. 

In this paper, therefore, we make a strengthened case for using the VHED functions introduced in \cite{greenleaf2018propagation} as inputs for the stroke classification problem. In particular, we continue the research line introduced in~\cite{agnelli2020classification} and extend the results to a scenario which is much more mathematically realistic and physically detailed than was employed in that earlier work:
(i) to simulate the virtual stroke patients, we consider phantoms which are more general than the unit disc, specifically, two-dimensional anatomical models resembling a slice of a human head with detailed, physically realistic electrophysiology and a statistically realistic conductivity distribution, and  
(ii) instead of using the unrealistic continuum model for EIT measurements, our new method uses realistic complete electrode model (CEM) measurements as its starting point.

Since the theory of VHED functions is based on measurements associated with the idealized continuum model, to compute the VHED functions from CEM measurements we first need to mimic continuous measurements from data originated using CEM. For this step, we follow the theory presented in~\cite{garde2021mimicking} and implement an algorithm that allows us to approximate unrealistic continuous measurements in two-dimensional EIT based on realistic electrode measurements.

This paper is organized as follows. 
In Section~\ref{sec:method}, we present the computation of the VHED functions from EIT data using the continuum model and then the extension in the case of the CEM.  Also, we explain our simulation of virtual patients, introduce the neural network architecture and describe the stroke classification problem.
Then, in Section~\ref{sec:resuls-discussion} the performance of the proposed methodology is evaluated using various case studies, and is compared against the case where raw EIT voltages are used as network inputs. Finally, conclusions are given in Section~\ref{sec:conclusions}.

\section{Methods}\label{sec:method}

In this section, we present how the VHED functions are computed from realistic EIT electrode measurements, the creation of the virtual patients with enhanced physiological detail, and the stroke classification method.

\subsection{From EIT data to VHED functions}

\subsubsection{EIT conductivity equation}
Let $\Om \subset \R^2$ be an open, bounded and simply connected domain (target), let $\sigma: \Om\to (0,\infty)$ be an essentially bounded measurable function (conductivity) satisfying $\sigma(x)\geq c>0$ for almost every $x\in \Om$, let $u \in H^1(\Om)$ be the voltage in the target, and let $\phi \in H^{1/2}(\dOm)$ be the voltage at the boundary of the target. Then the following conductivity equation holds \cite{somersalo1992existence}:
\begin{align}
\nabla \cdot \sigma \nabla u &= 0    \quad  \text{in  } \Om \label{eq:conductivity} \\
u &= \phi \quad \text{on } \dOm, \label{eq:dirichlet}
\end{align}
and the associated Dirichlet-to-Neumann (DN) map $\Lambda_{\sigma}$
is defined as
\begin{equation}\label{eq:DN-map}
\Lambda_{\sigma}: \phi  \longmapsto \sigma \frac{\partial u}{\partial \nu},
\end{equation}
where $\nu$ denotes the unit outer normal to the boundary. EIT images are generated by solving the inverse conductivity problem:
find $\sigma$ from the knowledge of $\Lambda_{\sigma}$ \cite{calderon_inverse_1980, astala2006calderon}.

\subsubsection{From continuum data to VHED profiles}\label{sec:CM-to-VHED}

There are many theoretical tools which have been developed to solve the conductivity equation~\eqref{eq:conductivity}, but here we focus on a special class of solutions called complex geometrical optics (CGO) solutions. These solutions are characterized by an exponential asymptotic behavior at infinity and provide the theoretical framework to introduce the VHED functions.
From now on, we assume that $\sigma=1$ outside of an open set $\Om_0$ satisfying $\overline{\Om_0} \subset \Om$. 

In~\cite{astala2006calderon} the authors present a way to construct the CGO solutions by reducing the conductivity equation~\eqref{eq:conductivity} to the following Beltrami equation:
\begin{equation}\label{eq:Beltrami}
	\dbar_z f_{\pm \mu} = \pm \mu\overline{\partial_z f_{\pm \mu}, }
\end{equation}
with Beltrami coefficient $ \mu := (1-\sigma)/(1+\sigma)$ and
the complex (Wirtinger) derivatives, $\partial_z=\tfrac 12(\partial_{x_1}-i\partial_{x_2})$
and  $\overline \partial_z=\tfrac 12(\partial_{x_1}+i\partial_{x_2})$.
The CGO solutions have the form
\[ f_{\pm \mu}(z,k) = e^{ikz}(1+\omega^{\pm}(z,k)),\] 
with  
\[ \omega^\pm(z,k)={\mathcal O}(\frac{1}{z}) \quad \mbox{as} \quad |z|\to \infty, \]
where $ikz = ik(x_1 + ix_2 )$, and where $z$ is considered a spatial
variable and $k$ a spectral parameter. Existence, uniqueness and properties of the CGO solutions are discussed in~\cite{astala2006calderon}.

Note that in order to compute $f_{\pm \mu}$, we need the coefficient $\mu$ which is computed from the conductivity $\sigma$. However, if we denote $W_{\pm}(z,k) = 1+\omega^{\pm}(z,k)$, by~\cite{astala2006boundary}, the following boundary
integral equation (BIE) holds:
\begin{equation}\label{eq:BIE}
	W_{\pm}(\,\cdot\,,k)|_{\partial\Om}+1 = ({\mathcal P}_{\pm }^k + {\mathcal P}_0)W_{\pm}(\,\cdot\,,k)|_{\partial\Om},
\end{equation}
where ${\mathcal P}_{\pm}^k$ and ${\mathcal P}_0$ are particular projection operators that can be  explicitly obtained from $\Lambda_{\sigma}$. Therefore, we can compute the traces of the CGO functions $\omega^{\pm}$ directly from  $\Lambda_{\sigma}$ by solving~\eqref{eq:BIE} and without the knowledge of $\mu$.

In\cite{greenleaf2018propagation}, a novel method connecting EIT and X-ray tomography
was introduced. This method converts information contained in the ideal continuum  DN map $\Lambda_{\sigma}$  into virtual X-ray profiles known as VHED functions, which contain specific geometric information about the conductivity.
In this methodology, the Beltrami equation~\eqref{eq:Beltrami} is treated as a scattering equation and the CGO function $\omega^{\pm}$ can be expanded as follows:
\[  \omega^{\pm}(z,k)= \sum_{n=1}^{\infty}  \omega^{\pm}_n(z,k).\]
The computation of the  VHED functions from the traces of the CGO functions $\omega^{\pm}$ consists of the following steps:

\begin{itemize}
	
	\item[(i)] 
	For any $\tau \in \R$ and $\varphi \in [0,2\pi]$, write the spectral variable of CGO solutions as $k=\tau e^{i \varphi}$. This is similar to polar coordinates, but not quite so as $\tau$ can be negative.
	
	\item [(ii)] Compute the partial Fourier transform in $\tau$, 
	\[
	\widehat{\omega}_n^{\pm}(z,t,\varphi):= \F_{\tau \to t}(\omega_n^{\pm}(z,\tau e^{i\varphi})).
	\]
	
	\item[(iii)] Compute the averaged operators $T_n^{\pm}\mu$ given by the complex contour integrals
	\begin{equation}\label{eq:T_op}
		T_n^{\pm}\mu(t,\varphi) := \frac{1}{2\pi i} \int_{\dOm} \widehat{\omega}_n^{\pm}(z,t,\varphi)  \, dz.
	\end{equation}
	
\end{itemize}

\begin{figure}
	\centering
	\includegraphics[width=4.8cm]{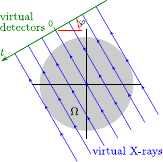}
	\caption{
	\small{
	Parallel VHED X-rays and virtual detector geometry.
    }
    }
	\label{fig:VHED_diredction}
\end{figure}

Note that in step (ii) we have introduced a nonphysical artificial (i.e., virtual) variable, $t$, known as the \emph{pseudo-time}. Fig.~\ref{fig:VHED_diredction} depicts the resulting virtual X-ray geometry. The angle $\varphi$ represents the direction perpendicular to the virtual X-rays (blue) while $t$ represents the position of a virtual X-ray detector, along the perpendicular direction (green).

\begin{figure}
	\centering
	\setlength{\unitlength}{1cm}
	\begin{picture}(10,6.35)
	
	\put(0.29,3.0){\includegraphics[scale=0.23,trim=60 0 70 0,clip]{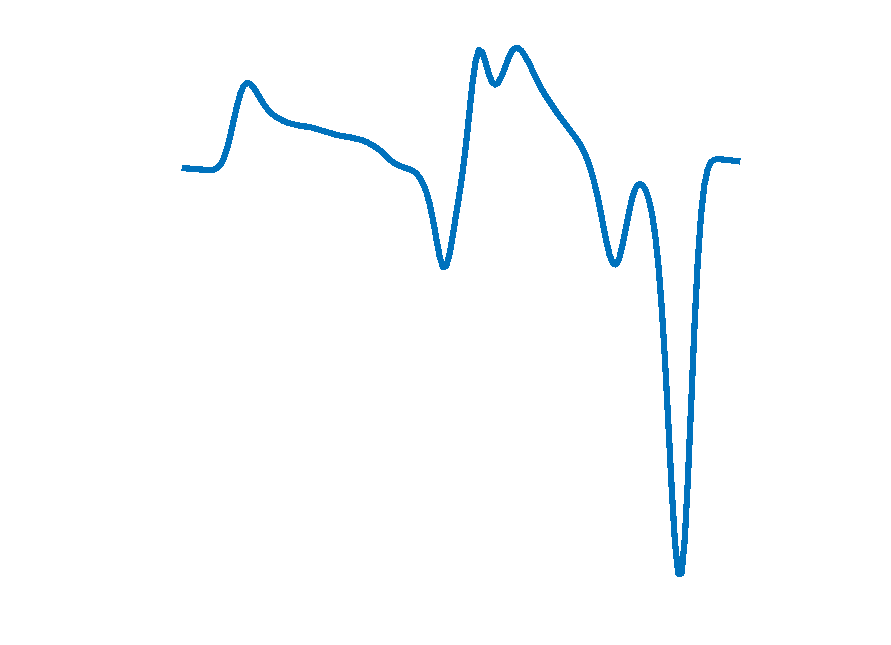}}
	\put(3.1,2.9){\includegraphics[scale=0.23,trim=52 0 52 0,clip]{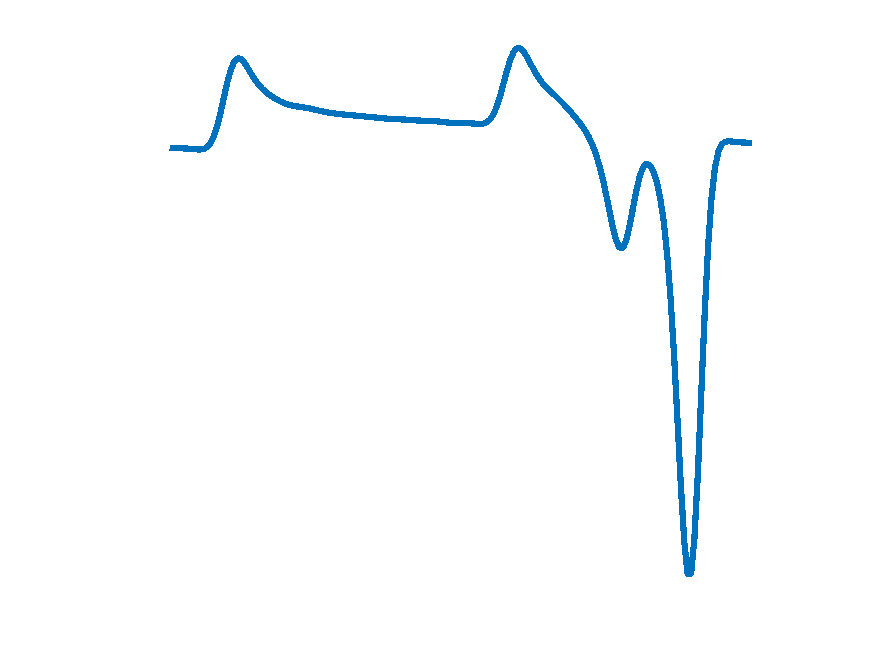}}
	\put(6.0,3.75){\includegraphics[scale=0.20, trim=25 0  40  0,clip]{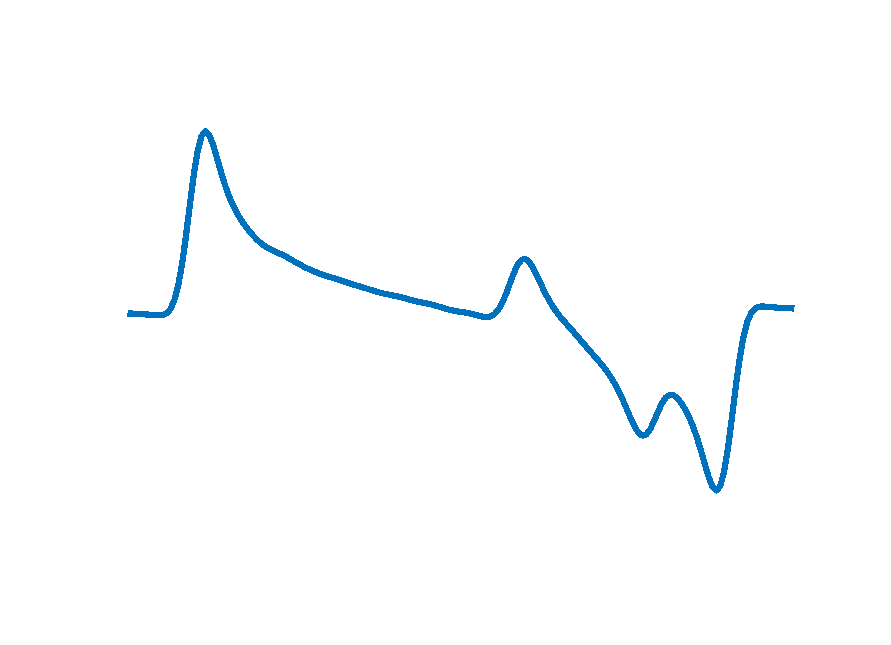}}
	\put(0.2,0.2){\includegraphics[scale=0.24,trim=52 0 52 0,clip]{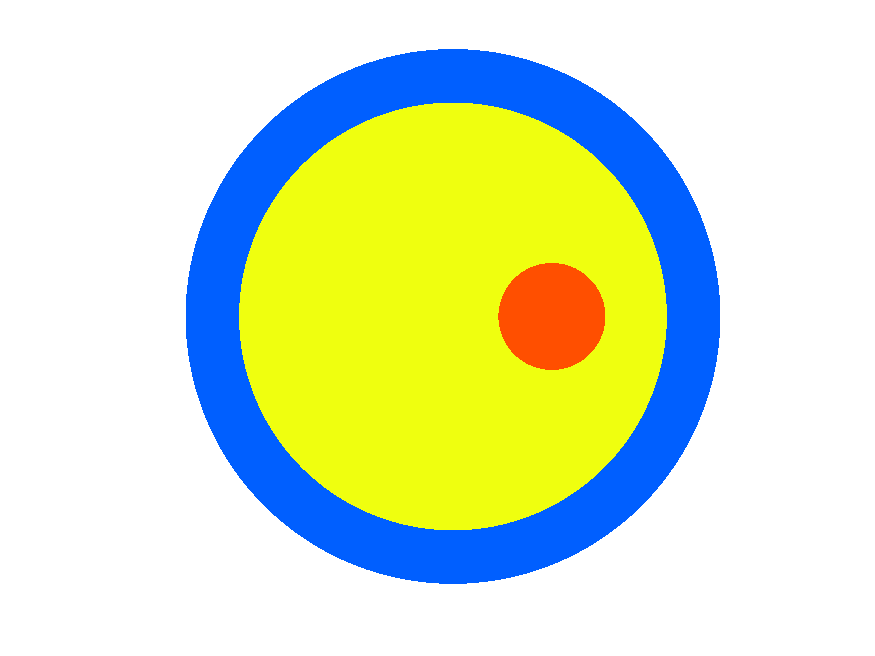}}
	\put(3.1,0.2){\includegraphics[scale=0.24,trim=52 0 52 0,clip]{images/Fig1_conductivity.png}}
	\put(6.0,0.2){\includegraphics[scale=0.24,trim=52 0 52 0,clip]{images/Fig1_conductivity.png}}
	
	
	\put(0.6,0.5){
		\begin{tikzpicture}
		\draw[dashed]  (.0,1.0)--(.0,3.5);
		\end{tikzpicture}
	}
	\put(1.65,0.5){
		\begin{tikzpicture}
		\draw[dashed]  (.0,1.0)--(.0,3.5);
		\end{tikzpicture}
	}
	\put(2.1,0.5){
		\begin{tikzpicture}
		\draw[dashed]  (.0,1.0)--(.0,3.5);
		\end{tikzpicture}
	}
	\put(2.33,0.5){
		\begin{tikzpicture}
		\draw[dashed]  (.0,1.0)--(.0,3.5);
		\end{tikzpicture}
	}

	\put(0.6,3.25){
		\begin{tikzpicture}
		\draw[dashed]  (.0,1.0)--(.0,3.6);
		\end{tikzpicture}
	}
	\put(1.65,3.25){
		\begin{tikzpicture}
		\draw[dashed]  (.0,1.0)--(.0,3.6);
		\end{tikzpicture}
	}
	\put(2.1,3.25){
		\begin{tikzpicture}
		\draw[dashed]  (.0,1.0)--(.0,3.6);
		\end{tikzpicture}
	}
	\put(2.33,3.25){
		\begin{tikzpicture}
		\draw[dashed]  (.0,1.0)--(.0,3.6);
		\end{tikzpicture}
	}
	
	\put(3.5,0.5){
		\begin{tikzpicture}
		\draw[dashed]  (.0,1.0)--(.0,3.5);
		\end{tikzpicture}
	}
	\put(4.55,0.5){
		\begin{tikzpicture}
		\draw[dashed]  (.0,1.0)--(.0,3.5);
		\end{tikzpicture}
	}
	\put(4.98,0.5){
		\begin{tikzpicture}
		\draw[dashed]  (.0,1.0)--(.0,3.5);
		\end{tikzpicture}
	}
	\put(5.25,0.5){
		\begin{tikzpicture}
		\draw[dashed]  (.0,1.0)--(.0,3.5);
		\end{tikzpicture}
	}
	
	\put(3.5,3.25){
		\begin{tikzpicture}
		\draw[dashed]  (.0,1.0)--(.0,3.6);
		\end{tikzpicture}
	}
	\put(4.55,3.25){
		\begin{tikzpicture}
		\draw[dashed]  (.0,1.0)--(.0,3.6);
		\end{tikzpicture}
	}
	\put(4.98,3.25){
		\begin{tikzpicture}
		\draw[dashed]  (.0,1.0)--(.0,3.6);
		\end{tikzpicture}
	}
	\put(5.25,3.25){
		\begin{tikzpicture}
		\draw[dashed]  (.0,1.0)--(.0,3.6);
		\end{tikzpicture}
	}
	
	\put(6.4,0.5){
		\begin{tikzpicture}
		\draw[dashed]  (.0,1.0)--(.0,3.5);
		\end{tikzpicture}
	}
	\put(7.47,0.5){
		\begin{tikzpicture}
		\draw[dashed]  (.0,1.0)--(.0,3.5);
		\end{tikzpicture}
	}
	\put(7.88,0.5){
		\begin{tikzpicture}
		\draw[dashed]  (.0,1.0)--(.0,3.5);
		\end{tikzpicture}
	}
	\put(8.15,0.5){
		\begin{tikzpicture}
		\draw[dashed]  (.0,1.0)--(.0,3.5);
		\end{tikzpicture}
	}
	
	\put(6.4,3.25){
		\begin{tikzpicture}
		\draw[dashed]  (.0,1.0)--(.0,3.6);
		\end{tikzpicture}
	}
	\put(7.47,3.25){
		\begin{tikzpicture}
		\draw[dashed]  (.0,1.0)--(.0,3.6);
		\end{tikzpicture}
	}
	\put(7.88,3.25){
		\begin{tikzpicture}
		\draw[dashed]  (.0,1.0)--(.0,3.6);
		\end{tikzpicture}
	}
	\put(8.15,3.25){
		\begin{tikzpicture}
		\draw[dashed]  (.0,1.0)--(.0,3.6);
		\end{tikzpicture}
	}
	
	\put(0.8,6.1){$\widehat{\omega}^+(1,t,\!0)$}
	\put(3.6,6.1){$\widehat{\omega}_{\mbox{\tiny odd}}^+(1,\!t,\!0)$}
	\put(6.5,6.1){$T_{\mbox{\tiny odd}}^{+}\mu (t,\!0)$}
	
	\put(1.0,4.1){\vector(1,1){0.4}}
	
	\put(1.35,0){(a)}
	\put(4.3,0){(b)}
	\put(7.15,0){(c)}
	
	\end{picture}    
	\caption{
		\small{ 
		VHED functions revealing geometric information about the conductivity.
		Jump singularities in the conductivity are reflected in the VHED profiles in a similar manner to parallel beam X-ray tomography.
		The angle $\varphi=0$ indicates the direction perpendicular to the virtual X-rays (dashed lines). 
		(a) Function $\widehat{\omega}^{+}(1,t,0)$ propagates singularities from the interior of the domain to its boundary. Note that an artifact is visible at $t=0$ (indicated by an arrow). 
		(b) Function $ \widehat{\omega}_{\mbox{\tiny odd}}^{+} := \left( \widehat{\omega}^+ - \widehat{\omega}^- \right) / 2$. Subtracting $\widehat{\omega}^{-}$ from $\widehat{\omega}^{+}$ eliminates the even terms and suppresses the artifact at $t=0$.
		Nevertheless, because of the election $z=1$, singularities closer to the right edge have a larger amplitude than the singularities closer to the left edge.
		(c) After applying the averaging operator, the VHED profile $T_{\mbox{\tiny odd}}^{+}\mu$ shows equal-size (but opposite sign) singularities at the leftmost and rightmost vertical lines.
	}
	}
	\label{fig:VHED_functions}
\end{figure}

This step permits that singularities from the interior of $\Om$ propagate to the boundary $\dOm$. Therefore, interior singularities of the unknown conductivity $\sigma$ can be detected using boundary measurements, see 
Fig.~\ref{fig:VHED_functions}(a).

Based on~\cite{greenleaf2018propagation}, the leading-order term $\omega^{\pm}_1$ contains relevant information that allows detection of edges and singularities of the coefficient $\mu$, and thus of the conductivity $\sigma$.
In fact, in~\cite{greenleaf2018propagation} it is shown that the following filtered back-projection formula holds:
\begin{equation*}\label{eq:FBP}
	(- \Delta)^{-1/2} (T_1^{\pm})^{\ast}T_1^{\pm} \mu = \mu
\end{equation*}
and 
\[ 
T_1^{\pm} \mu(t,\varphi) = -2\pi i e^{i\varphi} \big(\frac{\partial}{\partial s}R\mu_{\pm}\big)(\frac12 t,e^{i\varphi}),
\]
where $R$ is the standard Radon transform. 

Unfortunately, we can only compute the boundary trace of the full scattering series $\omega^{\pm}$  from the DN map. 
As we mention above, the first-order term $\omega^{\pm}_1$ constitutes the basis for an accurate reconstruction of the unknown conductivity $\sigma$, while the higher-order terms $\omega_n^{\pm}$ are related to scattering effects and these terms produce artifacts in the numerical simulations, see Fig.~\ref{fig:VHED_functions}(a).
However, by a parity symmetry property described in~\cite{greenleaf2018propagation}, subtracting $\widehat{\omega}^{-}$ from $\widehat{\omega}^{+}$ eliminates the even terms, $\omega_{2n}^{+}$, so that their singularities can be suppressed, see Fig.~\ref{fig:VHED_functions}(b). In this study we will work with the following VHED function:
\[ T_{\mbox{\tiny odd}}\mu(t,\varphi) := \frac{T^{+}\mu(t,\varphi) - T^{-}\mu(t,\varphi)}{2}, \]
where $T^{\pm}\mu$ denotes the average operator~\eqref{eq:T_op} applied to 
$\widehat{\omega}^{\pm}$ respectively, see Fig.~\ref{fig:VHED_functions}(c).

\subsubsection{From electrode data to VHED profiles}\label{subsec:CEM-to-VHED}
Note that in order to compute the VHED functions, we need the DN map $\Lambda_{\sigma}$ associated with the continuum model~\eqref{eq:conductivity}-\eqref{eq:dirichlet}. However, in practical applications, a finite number of electrodes are attached at the boundary of the body. The resulting EIT measurement data comprise pairs of voltage and current measurements recorded at these electrodes. 
The \emph{complete electrode model} (CEM)~\cite{somersalo1992existence} is considered the most realistic mathematical representation for simulating electrode measurements.

Let $e_m \subset \dOm$, $m=1,\dots,M$ be disjoint open and connected sets modeling the electrode patches. The CEM is given by the conductivity equation~\eqref{eq:conductivity} and the following boundary conditions:
\begin{align*}
z_m  \, \sigma \frac{\dd u}{\dd \nu}  + u  &= U_m, \quad \,  \text{on  } e_{m}, \;  m=1,\ldots, M   \\
\sigma \frac{\dd u}{\dd \nu} &= 0, \qquad  \text{on }  \dOm \setminus \cup_{m = 1}^{M}e_m  \\
\int_{e_m} \sigma \frac{\dd u}{\dd \nu} \, ds &= I_m , \quad \, \,  m=1,\dots,M,
\end{align*}
where $z_m$ denotes the contact impedance at electrode $e_m$,
$U_m$ represents the constant electric potential on electrode $e_m$, $I_m$ is the net current through the electrodes, and the normal current density outside the electrodes is zero.

Within the CEM framework, the \emph{absolute} measurements are modeled by the electrode current-to-voltage map
$\cR_\sigma^{\text{CEM}}:  \R^M_{\diamond} \to  \R^M_{\diamond}$ given by 
\begin{equation}\label{eq:ND-map-CEM}
	\cR_{\sigma}^{\text{CEM}}: I  \longmapsto  U,
\end{equation}
where $\R^M_{\diamond}$ denotes the mean free subspace of $\R^M$.
Imposing the current conservation law
$\sum_{m=1}^{M} I_m = 0$
and fixing the ground level of the potentials
$\sum_{m=1}^{M} U_m = 0$ 
guarantees the existence and uniqueness of the solution $(u,U)$ of the CEM, where $u \in H^1(\Om)$ and $U =(U_1, \ldots, U_M )^T \in \R^M_{\diamond}$, for details see~\cite{somersalo1992existence}.

In~\cite{garde2021mimicking} a method for approximating relative continuum EIT data based on measurement data from the CEM is presented. Given a bounded simply connected domain $\Om \subset \R^2$ and a set of continuum current patterns $\{\phi_k\}_{k=1}^{M-1}$
that one would like to apply along the boundary $\dOm$, the proposed method tells us where the electrodes, $e_m$, should be located, based on the shape of $\Om$, and what the net electrode currents  $\{ I^{(k)} \}_{k=1}^{M-1}$
should be, depending on the continuum current patterns. 

The method is based on accurate interpolation and quadrature rules on the boundary of the unit disc, $D \subset \R^2$, for certain Sobolev spaces. Then, with the help of the Riemann mapping theorem, these results are extended to more general domains.
Next, we present the algorithm that allows us to compute a matrix approximation $\mathbf{L}_{\sigma}$ of the continuum DN map $\Lambda_{\sigma}$ given by~\eqref{eq:DN-map} but using boundary measurements from the CEM model instead of the continuum model. 

It is worth mentioning that in~\cite{garde2021mimicking}, instead of considering the Dirichlet condition~\eqref{eq:dirichlet}, the Neumann condition 
\[
\frac{\dd u}{\dd \nu} = \phi \quad \text{on } \dOm \label{eq:neumann}
\]
is considered. 
In this case, the boundary measurements are modeled by the Neumann-to-Dirichlet (ND) map
\begin{equation*}
\cR_{\sigma}: \phi  \longmapsto  u_{|_{\dOm}},
\end{equation*} 
and relative continuum data associated to the map $\cR_{\sigma}-\cR_{1}$ is approximated based on CEM measurements, where $\cR_{1}$ is the ND map with $\sigma=1$ everywhere. However, once a matrix approximation $\mathbf{R}_{\sigma}$ of the ND map is available, the computation of a matrix approximation $\mathbf{L}_{\sigma}$ of the DN map $\Lambda_{\sigma}$ is done using the inverse $\mathbf{R}_{\sigma}^{-1}$.

The first step of the constructive method presented in~\cite{garde2021mimicking} consists of locating the electrodes $e_m$ in precisely defined places.   
Let $\displaystyle{\theta_m = m \frac{2 \pi}{M}}$ and $x_m = \text{e}^{i\theta_m}$, for $m=1,\dots,M$. The points $x_m$ can be interpreted as the electrodes' midpoints on the boundary of the unit disc $D$. 
Let $\Psi: D \to \Omega$ be a bijective conformal map whose restriction to $\dd D$ defines a $C^{\infty}$-diffeomorphism between $\dd D$ and $\dOm$. 
Denote by $y_m := \Psi(x_m)$, for $m=1,\dots,M$. The points $y_m$ represent the electrodes' midpoints on the boundary of the domain $\Omega$. 
Additionally, we assume that we are given an arc length parametrization $\gamma:(0,L] \to \partial \Omega$, where $L$ is the total arc length of $\partial \Omega$,  and we denote by $s_m$ the points such that $\gamma(s_m) = y_m= \Psi(x_m),$ for $m=1,\dots,M$. 

Once we have defined the electrode locations, the second step consists of defining the set of input currents given the continuum current patterns. In this study we consider, as the continuum patterns, the following set of trigonometric basis functions 
\begin{equation} \label{eq:trig-cp}
\phi_k(s)= 
\begin{cases} 
\displaystyle \sqrt{\frac{2}{L}} \cos\left(\frac{(k+1)}{2}\frac{2\pi s}{L} \right) &\quad \mbox{odd } k, \\
\displaystyle \sqrt{\frac{2}{L}} \sin\left(\frac{k}{2} \frac{2\pi s}{L} \right) &\quad \mbox{even } k,
\end{cases}
\end{equation}
for $k=1,\dots,M-1$. The corresponding electrode current pattern for the CEM is defined as
\begin{equation*}\label{eq:I_cp}
I^{(k)} \! = \!  \frac{2\pi}{M} \!
\left( \!
\begin{bmatrix}
\phi_k(s_{1}) \, |\Psi'(x_{1})| \\
\vdots\\
\phi_k(s_{M}) \, |\Psi'(x_{M})| 
\end{bmatrix} \! - \!  \frac{1}{M} \sum_{m=1}^{M} \! \phi_k(s_m) \, |\Psi'(x_{m})| \,  \mbox{\bf{1}} \!
\right)
\end{equation*}
for $k=1,\dots,M-1$ and where $\mbox{\bf{1}}=[1,\ldots,1]^T \in \R^M.$

Finally, the process to construct a matrix approximation $\mathbf{L}_{\sigma}$ of the DN map $\Lambda_{\sigma}$ consists of: 
\begin{enumerate}
	
	\item  For $k=1,\ldots,M-1$, take $I^{(k)}$ as input current, perform the measurements of CEM to get $ U^{(k)}= \left(\cR_{\sigma}^{\text{CEM}} -\cR_{1}^{\text{CEM}} \right) I^{(k)} \in \R^{M}_\diamond$ and then form the matrix 
	\[ \mathbf{U} = [ U^{(1)} \, U^{(2)} \ldots U^{(M-1)}] . \]

	\item  Construct the matrix
	\[ \mathbf{I} = [ I^{(1)} \, I^{(2)} \ldots I^{(M-1)}] .\]
	
	\item Compute a matrix approximation of the relative electrode current-to-voltage map as
	\[ \left( \mathbf{R}_{\sigma}^{\text{CEM}} - \mathbf{R}_{1}^{\text{CEM}} \right) = \mathbf{U}^T \mathbf{I}. \]

	\item Compute a matrix approximation $\mathbf{R}_{1}$ of $\cR_1$ in the basis~\eqref{eq:trig-cp}.

	\item Compute a matrix approximation $\mathbf{\widetilde{L}}_{\sigma}$ of the DN map $\Lambda_{\sigma}$ in the basis~\eqref{eq:trig-cp} as 
	\begin{equation*}
	\mathbf{\widetilde{L}}_{\sigma} = ( \mathbf{U}^T \mathbf{I} + \mathbf{R}_{1}) ^{-1}.
	\end{equation*}
\end{enumerate}
However, the DN map $\Lambda_{\sigma}$ acts also on constant functions.
Therefore, we extend the matrix $\widetilde{\mathbf{L}}_{\sigma}$  and define
\begin{equation*}\label{eq:L-sigma}
\mathbf{L}_{\sigma} = 
\begin{bmatrix}
0 & 0 \\
0 & \widetilde{\mathbf{L}}_{\sigma} 
\end{bmatrix},
\end{equation*}
obtaining a matrix approximation of the DN map $\Lambda_{\sigma}$ acting on the basis~\eqref{eq:trig-cp} augmented by the constant basis function $\phi_0(s) = |\dOm |^{-1/2}$.

\begin{figure}
	\centering
	\setlength{\unitlength}{1cm}
	\begin{picture}(10,1.5)
	\put(1,0){\includegraphics[scale=0.32,trim=25 0 40 0,clip]{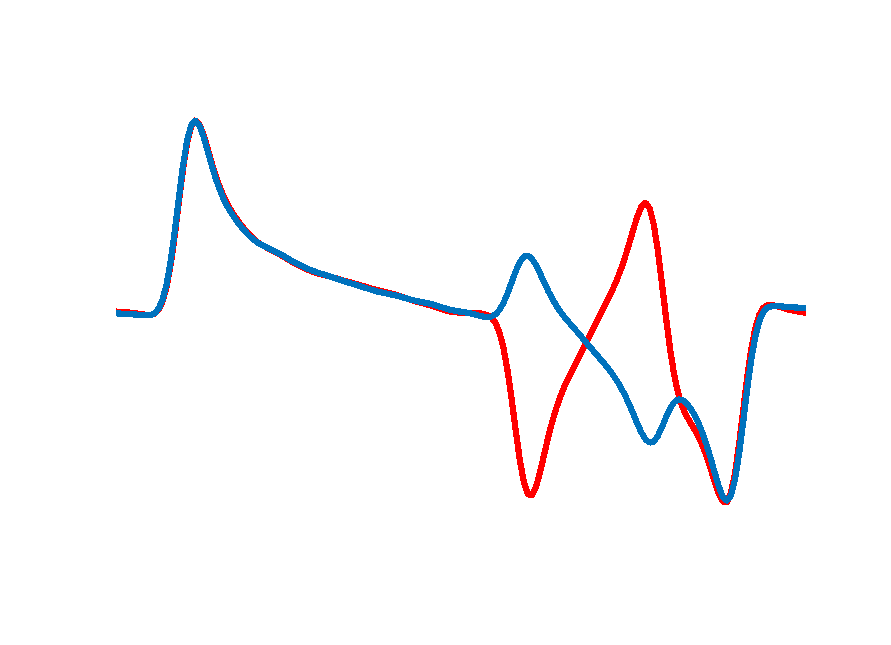}}
	\put(5.25,-0.25){\includegraphics[scale=0.32,trim=25 0 40 0,clip]{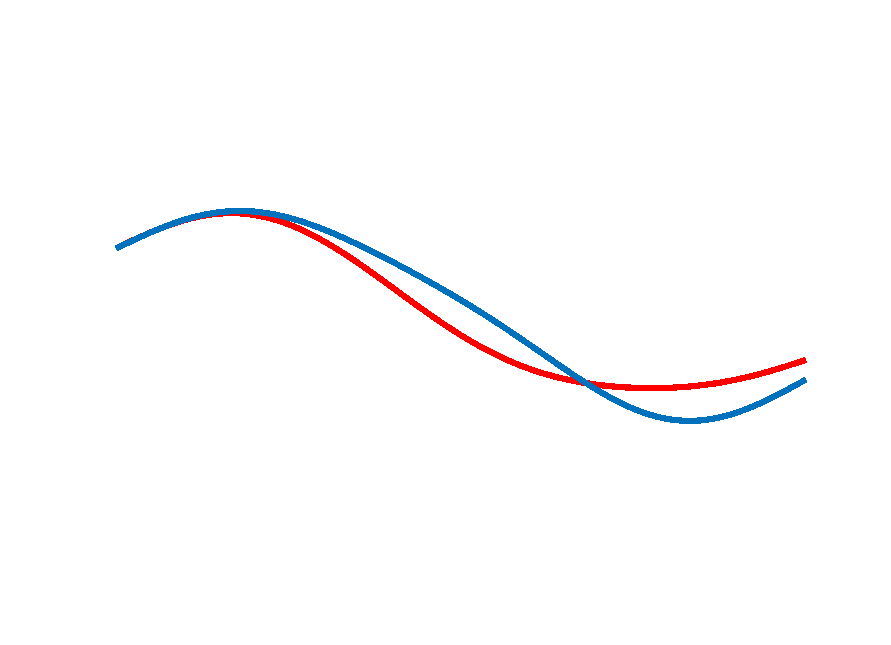}}
	
	\put(3.0,0){(a)}
	\put(7.5,0){(b)}
	
	\end{picture}    
	\caption{
	\small{	
		Effect of the Fourier windowing process applied to the VHED functions. The blue curve shows $T_{\mbox{\tiny odd}}\mu(t,0)$ corresponding to conductivity in Fig.~\ref{fig:VHED_functions} where the small disc inclusion has a conductivity value higher than background, while the red curve shows $T_{\mbox{\tiny odd}}\mu(t,0)$ corresponding to the same conductivity but with the small disc inclusion having lower conductivity than background. (a) Unrealistic VHED profiles with cutoff frequency $\tau_{m} =50$, which would only be possible in a noise-free scenario. (b) Realistic VHED profiles considering a cutoff frequency $\tau_{m}=5$, as would be appropriate in a real-world scenario with noise. 
	}
	}
	\label{fig:VHED_blurred}
\end{figure}

Once we have a matrix approximation $\mathbf{L}_{\sigma}$ of the DN map $\Lambda_{\sigma}$ we can use this to compute the traces of the CGO solutions by numerically solving the BIE~\eqref{eq:BIE}.
In practical applications, the presence of noise in the DN matrix leads to inaccuracies in the solution of~\eqref{eq:BIE}. From simulations, one can observe that errors in the values of the CGO solutions $\omega^{\pm}(z,\tau e^{i \varphi})$ for $z \in \dOm$ grow with $|\tau|$. Therefore, in order to mitigate inaccuracies when computing the VHED functions, we restrict the integration to the interval $[-\tau_{m},\tau_{m}]$ ($m$ for maximum value) and consider the following windowed Fourier transform
\begin{equation*}\label{eq:breve_w_1}
	\breve{\omega}^{\pm}(z,t,\varphi) := \int_{-\tau_{m}}^{\tau_{m}}  W_{a_{\tau_{m}}}(\tau) e^{-i \tau t } \omega^{\pm}(z,\tau e^{i\varphi}) \, d\tau,
\end{equation*}
where $a_{\tau_{m}} \! > \! 0$ is a parameter such that $W_{a_{\tau_{m}}}(\tau): = e^{ -a_{\tau_{m}} \tau^2 }$ is very small outside the interval $[-\tau_{m},\tau_{m}]$. 
Note that this Fourier windowing process in the frequency domain is equivalent to convolution in the pseudo-time domain, that is:
\[ \breve{\omega}^{\pm}(z,t,\phi) = \widehat{W}_{a_{\tau_{m}}}(t) *  \widehat{\omega}^{\pm}(z,t,\phi).\] 
Hence, due to the inevitable noise in the measurements we are compelled to work with blurred versions of the VHED functions see Fig.~\ref{fig:VHED_blurred}.
It is clear from Fig.~\ref{fig:VHED_blurred} that if we could use unrealistic noise-free data, the classification of strokes would be an easy task, and could be performed without the aid of machine learning. However, when using realistic  noisy data this task becomes much more complex, and so we rely upon a neural network to accurately perform classification.

\subsection{Simulation of stroke-EIT measurements}\label{subsec:head-model}
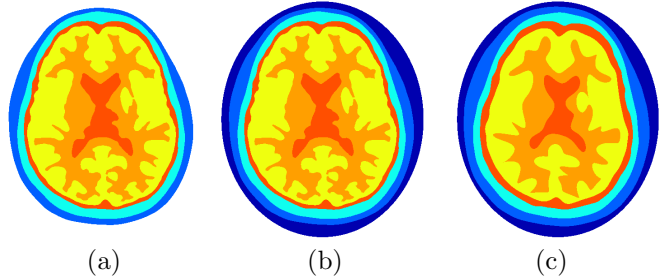
\begin{figure}
	\centering
	\setlength{\unitlength}{1cm}
	\begin{picture}(10,2.0)
	\put(0.3,0){\includegraphics[scale=0.35,trim=100 0 90 0,clip]{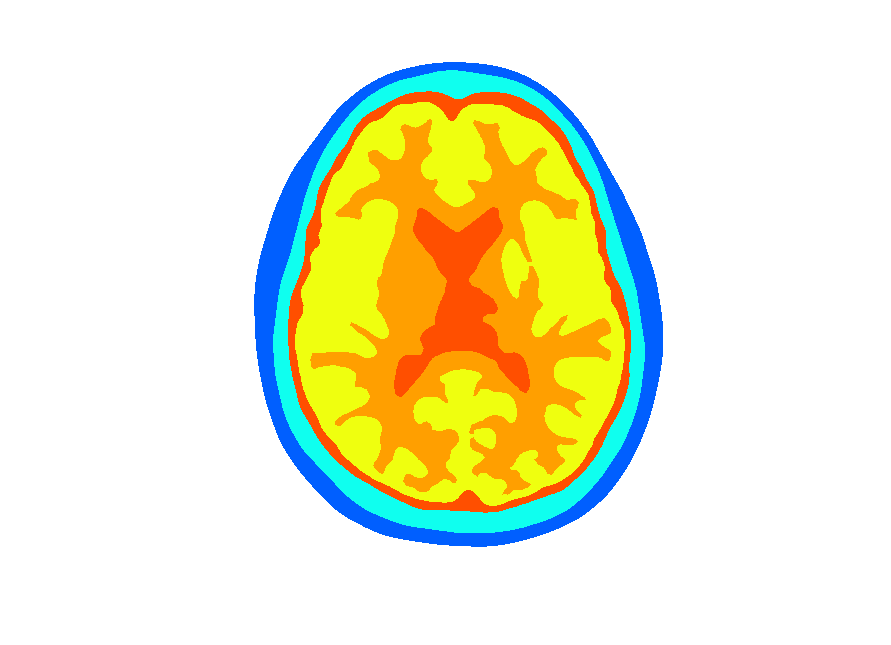}}
	\put(3.2,0){\includegraphics[scale=0.35,trim=100 0 90 0,clip]{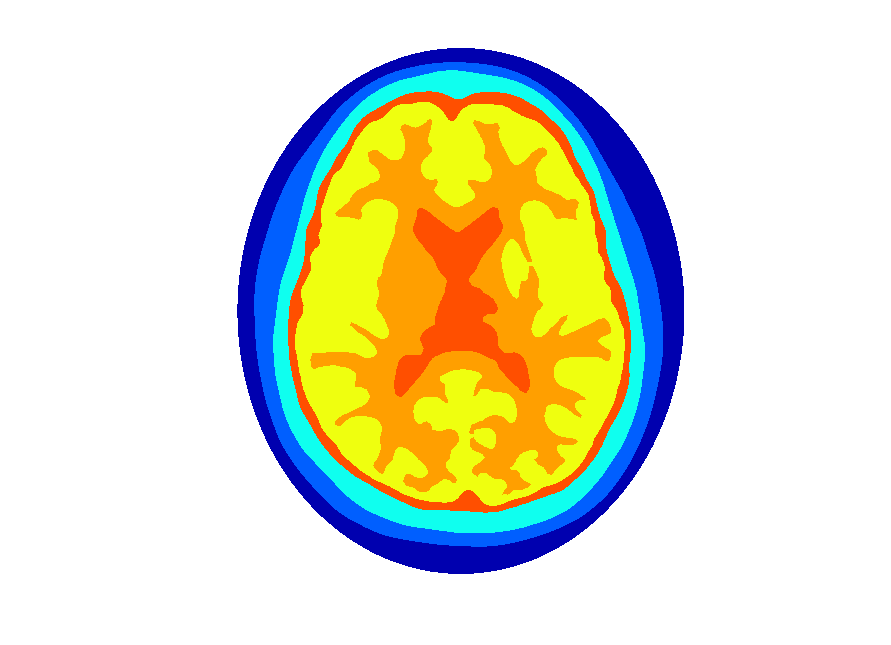}}
	\put(6.3,0){\includegraphics[scale=0.35,trim=100 0 90 0,clip]{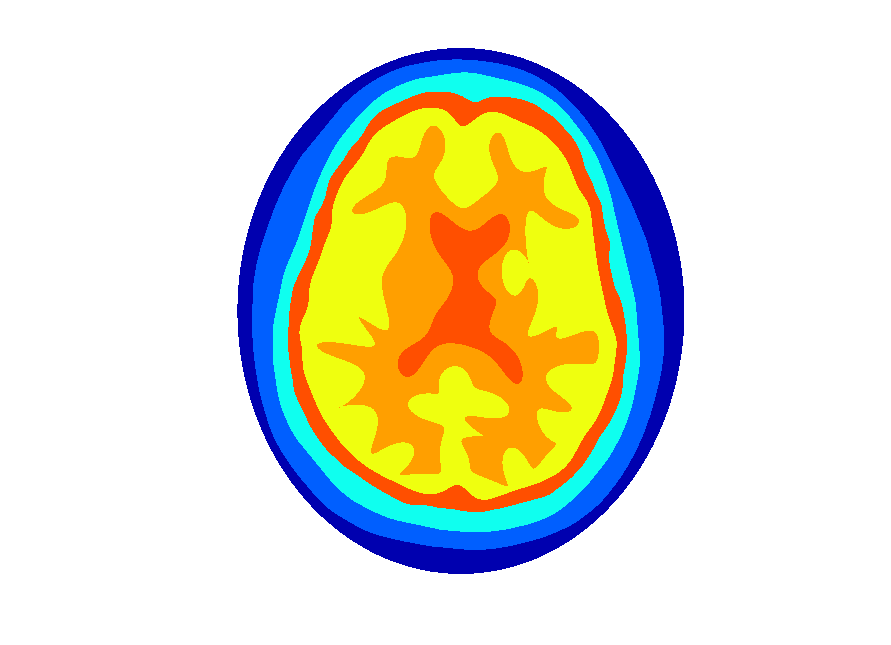}}
	
	\put(1.6,0){(a)}
	\put(4.5,0){(b)}
	\put(7.5,0){(c)}
	
	\end{picture}   
	\caption{
	\small{
	Two-dimensional model resembling a slice of a human head. (a) Reference model containing five different compartments of importance for electrophysiology of the human head: grey matter, white matter, cerebrospinal fluid, skull and other soft tissues representing the scalp. (b) Reference model included in a elliptic shape. The dark blue region can be interpreted as a fixed helmet of constant conductivity with electrodes attached to its exterior. (c) The head anatomy inside the helmet of one simulated patient.  
	}
	}
	\label{fig:head_models}
\end{figure}

A statistical anatomical atlas of the electrical properties of the upper part of the human head was introduced in~\cite{moura2021anatomical}.
Based on this work, we considered a physiological detailed two-dimensional  model resembling a slice of a human head. The model contains five different compartments of importance for electrophysiology of the human head: grey matter, white matter, cerebrospinal fluid (CSF), skull and other soft tissues representing the scalp, see 
Fig.~\ref{fig:head_models}(a). 

Recall that in order to compute the CGO solutions in an arbitrary domain $\Om \subset \mathbb{R}^2$ using CEM measurements, the first step is to find a bijective conformal map $\Psi$ from the unit disc $D$ to $\Om$. If we took $\Om$ to extend to the outer layer of the patient's skin, in practice this means that for each patient we would need to find a different conformal map, and this could be a difficult task. As a simple and physically realistic solution, therefore, we assume that a fixed elliptic belt or helmet with electrodes is placed on the head of the patient. In this way the boundary of the domain $\Om$ is fixed, only the inner compartments of the head model differ from patient to patient, see Fig.~\ref{fig:head_models}(b). Moreover, we selected an average head geometry as a reference  model, and in order to simulate different patients, we randomly modified the shape, width and location of the various components, see Fig.~\ref{fig:head_models}(c). Therefore, each sample has a unique head anatomy included in a fixed elliptic region with semi-axes $0.85$ and $1$. 

Additionally, in each of these head models the conductivity value assigned to each component of the head is randomly sampled from a Gaussian distribution with mean and standard deviation in the ranges given in Table~\ref{table:conductivities}. 
These realistic conductivity values were taken from~\cite{moura2021anatomical}. Note that all conductivity values are normalized by the conductivity of the scalp (this is needed for the computation of CGO solutions). The only region with constant conductivity in all samples is the blue region between the scalp and the boundary of the elliptic domain, which may be interpreted as a fixed helmet of constant conductivity.

{\small
	\ctable[
	mincapwidth = 8cm,
	caption = Normalized conductivity distribution of the different regions of the simulated two-dimensional head model.,
	label = table:conductivities,
	pos = h
	]{ll}{}{
		\FL
		Scalp              &  $\sim \N(1,0.0333)$       \NN
		Skull              &  $\sim \N(0.0625,0.0021)$  \NN
		Cerebrospinal fluid &  $\sim \N(6.25,0.2083)$    \NN
		Grey matter        &  $\sim \N(0.3063,0.0102)$	\NN
		White matter       &  $\sim \N(0.1938,0.0065)$	\NN
		Ischemic stroke    &  $\sim \N(0.0938, 0.0031)$	\NN
		Hemorrhagic stroke &  $\sim \N(2.1875, 0.0729)$	\LL
}}

It is worth noting that the conductivity of the skull is very low and therefore is highly resistive to the electric currents. On the other hand, the CSF is highly conductive and acts as a shunt path for the current. These two regions form a double shield against the passage of electric current through the brain; any current that crosses the skull encounters an easy path via the CSF. This important effect is taken into consideration in this work. One of the goals of this work, therefore, was to demonstrate that VHED functions computed from EIT measurements can be used for stroke classification even in the presence of these physically realistic structures, which are frequently omitted from models presented in other works, due to the challenges they create.  

Finally, stroke inclusions consisting of discs or ellipses of various eccentricities and orientations are considered, see Fig.~\ref{fig:stroke_models}. Exact parameters for the shapes and sizes of these inclusions, which differed for training and testing sets, are described in detail later in Section~\ref{sec:training-testing}. The conductivity values for ischemic (low-conductivity) vs. hemorrhagic (high-conductivity) strokes were drawn from the distributions shown in Table~\ref{table:conductivities}. We assume that stroke inclusions are always on the right-hand side of the head model. That is because, in general, strokes are restricted to one half of the brain only and the symptoms clearly indicate which half is affected. Thus, there is no loss of generality in locating all the simulated strokes on the right-hand side of the brain, as the same process could be done for strokes located on the left-hand side of the brain.

To generate the synthetic data needed for the classification problem we consider a large set of conductivity profiles with the properties described above, and for each sample we solve the CEM model, using the finite element method \cite{Brenner2008}, to obtain the corresponding EIT measurements. To recreate a more realistic scenario we add noise to the computed EIT data. That is, for each sample we collect all the electric potentials $U^{(k)}$, $k=1,\ldots,M-1$ and form the single real-valued column vector 
\[ V = [ U^{(1)}; \, U^{(2)}; \, \ldots \, ;U^{(M-1)}] \in \mathbb{R}^{M(M-1)}.\]
In practice, modern EIT systems tend to have $2^n$ voltage input channels, with some systems capable of handling measurements on up to 128 or 256 electrodes \cite{Xu2007, Aristovich2014, Cherepenin2002}. In this work, we followed the convention in \cite{garde2021mimicking} and used an odd number of electrodes, $M=65$, although this number is chosen for convenience and the calculations could easily be formulated to work with arbitrary $M$.

After obtaining voltages, we add zero-mean random Gaussian noise to the vector $V$ to obtain a noisy vector $V^{\delta}$, with relative noise given by
\[\frac{\|V-V^\delta\|_2}{\|V\|_2} = \delta,\]
where $\|\cdot\|_2$ denotes the $2$-norm. Following the techniques suggested in previous works, this noisy vector may be used directly as an input for the classification problem. In this present work, this method serves as a basis of comparison for our alternative strategy, wherein the voltages are further processed to obtain VHED functions. In this latter case, we reshape $V^\delta$ into a matrix and compute the noisy VHED profiles following the steps previously described in Section~\ref{subsec:CEM-to-VHED}. In this paper we consider $\delta=10^{-3}$ and $\delta=10^{-2}$, and, as appropriate for these levels of noise, we take a cutoff frequency $\tau_{m}=5$ in the computation of the VHED profiles. 

Computation of the VHED profiles is not computationally intensive, and so does not add significant time to the stroke classification process as compared to  the use of unprocessed voltages. In this work for example, using code which had not been optimized for efficiency, computation of VHED profiles from one set of EIT voltage data took only 28 seconds on an simple laptop with an Intel i7 processor and 16 GB of RAM, with the calculation of the CGO solutions being the most time-consuming part. However, the code could be optimized and run on a higher-performance machine to greatly reduce computation time, making the calculation and subsequent use of VHED functions a feasible diagnostic strategy for real-time applications.

\subsection{Stroke classification using neural networks}

The goal of this study is to explore the ability of machine learning to classify strokes as ischemic or hemorrhagic using EIT data simulated from mathematically realistic and physically detailed models.
As mentioned in Section~\ref{sec:intro}, most machine learning-based approaches proposed for the classification of strokes have considered learning directly from measured EIT voltage data~\cite{macdermott2018brain,McDermott2020,shi2022RCNN,candiani2022NN,Culpepper2024}. Here, we compare that method to the techniques which were introduced in~\cite{agnelli2020classification}, wherein noise-robust VHED functions computed from highly simplistic and mathematically ideal phantoms were used as network inputs.  Next, we describe the neural networks we employed to compare these two different paradigms, as well as the details of the training and testing datasets.

\subsubsection{Fully connected neural network (FCNN)}
In this work we consider a simple FCNN whose architecture consists of one input layer, one hidden layer and one output layer. It should be noted that because we were able to obtain excellent and robust results using this simple network, we did not pursue the use of more sophisticated machine learning techniques for this study.

The size of the input layer depends on the type of input data. 
If raw voltage data is considered as input, then the potential measurements are stacked into a single real-valued column vector 
$V \in \mathbb{R}^{M(M-1)}$
which is used as the input. As mentioned in Section~\ref{subsec:head-model} in this work we consider $M=65$ and therefore each input $V \in \mathbb{R}^{4160}$.
If a VHED function $T_{\mbox{\tiny odd}}\mu(t,\varphi)$ is instead considered as input data, this function is computed at angles $\varphi_i = i \frac{\pi}{8}$, for $i=4,5,6,12,13,14$, and the pseudo-time $t \in [-3,3]$ is discretized over $256$ points, so that for each angle $\varphi_i$ we have $T_{\mbox{\tiny odd}}\mu(t,\varphi_i) \in \mathbb{C}^{256}$. Then, the real and imaginary parts of the VHED profile at angles $\varphi_i$ are concatenated to form a single real-valued column vector $ T_{\mbox{\tiny odd}}\mu \in \R^{3072}$ that will be used as input. In order to unify notation we denote both inputs as $x \in \R^{n_1}$.

The hidden layer has 30 neurons. This number was selected based on the results obtained in preliminary tests. The output layer consists of one neuron. The output, denoted by $y$, is a scalar  that takes values in $\{0,1\}$, where $y=0$ represents an ischemic inclusion and $y=1$ represents a hemorrhagic inclusion.

From a mathematical point of view, the FCNN can be represented as function $\F_{\theta}:\mathbb{R}^{n_1}\to \mathbb{R}$ given by 
\[
x \longmapsto \F_{\theta}(x) := f(W_{\theta}^2f(W_{\theta}^1x +b_{\theta}^{1})+b_{\theta}^2),
\]
where $W^{1}_{\theta} \in \mathbb{R}^{30 \times n_1}$ and $b^{1}_{\theta} \in \mathbb{R}^{30}$ represent the weight matrix and the vector of biases at the hidden layer, respectively, $W_{\theta}^{2} \in \mathbb{R}^{1 \times 30}$ and $b_{\theta}^{2} \in \mathbb{R}$ represent the weight and the bias of the output layer, and  $f$ is the sigmoid activation function $f(t)=(1+e^{-t})^{-1}$, applied element-wise. Using a more compact notation, the FCNN can be represented as
\[ y=\F_{\theta}(x), \]
where $\theta \in \mathbb{R}^{30n_1+30+30+1}$ is a vector containing all the network parameters. In order to find the optimal parameters of the network we consider the minimization of the binary cross entropy loss function,
\begin{equation*}\label{bi_cross_entropy}
{\cal L} (\theta) = \sum_{j}-y_j \text{log}(y^p_j) - (1-y_j) \text{log}(1-y^p_j),
\end{equation*}
where $y_j$ is the true value of the outcome, $y^p_j= \mathcal{F}_{\theta}(x_j)$ is the predicted value of the outcome, and the sum is over all samples in the training dataset. The scaled conjugate gradient algorithm was considered for minimization of the loss function.

\subsubsection{Training and testing datasets}\label{sec:training-testing}

In this study, we consider that an elliptic belt or helmet with electrodes is placed on the head of the patient and we assume that the measurements come from a two-dimensional cross-section of the human head. As we described in Section~\ref{subsec:head-model}, a cross-section containing five different compartments of importance for electrophysiology is considered as the head model, and an inclusion having low or high conductivity simulates an ischemic or hemorrhagic stroke respectively.

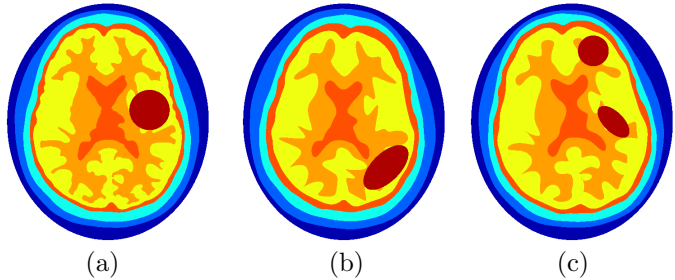
\begin{figure}
	\centering
	\setlength{\unitlength}{1cm}
	\begin{picture}(10,2.75)
	\put(0.5,0){\includegraphics[scale=0.35,trim=100 0 90 0,clip]{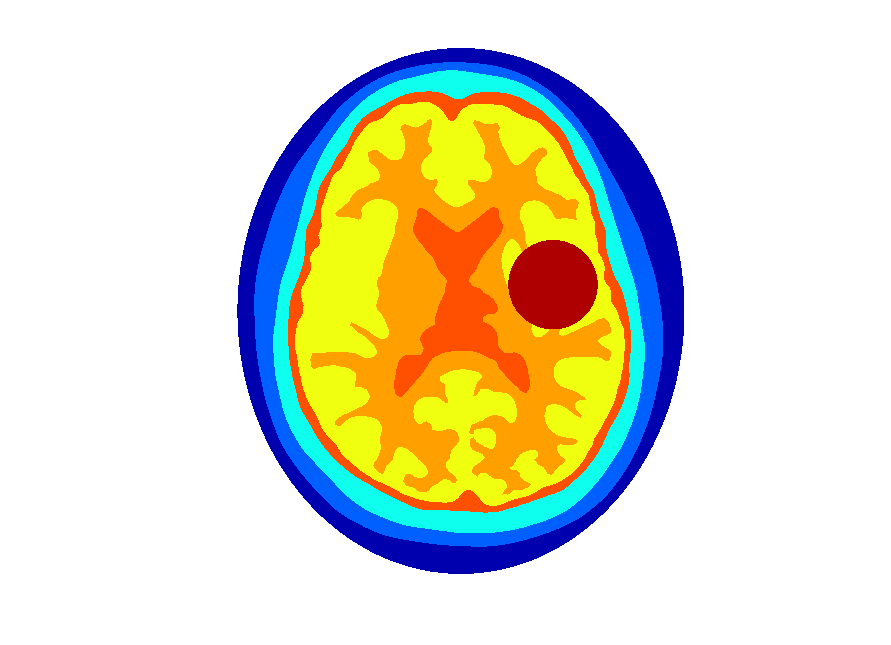}}
	\put(3.6,0){\includegraphics[scale=0.35,trim=100 0 90 0,clip]{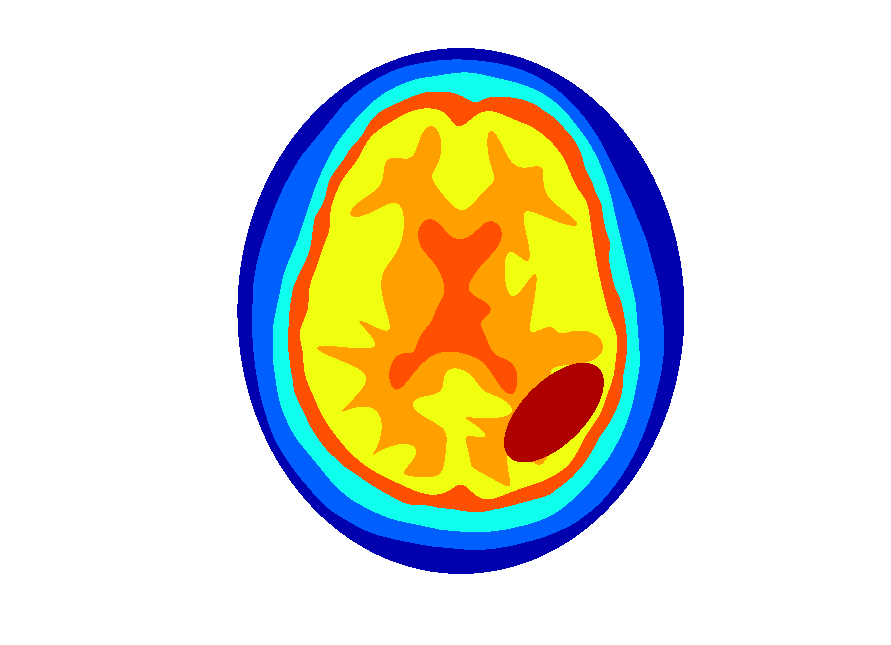}}
	\put(6.6,0){\includegraphics[scale=0.35,trim=100 0 90 0,clip]{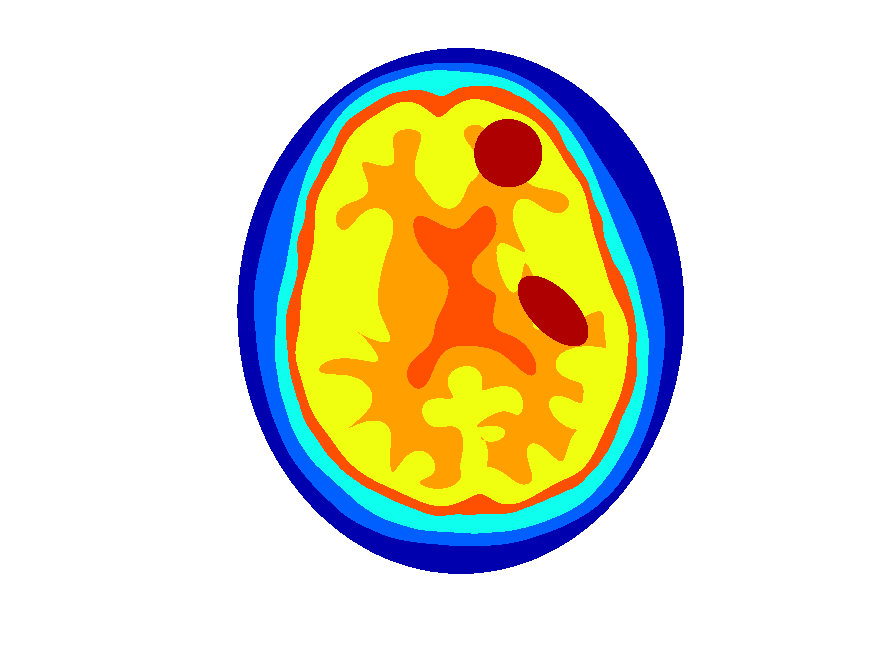}}
	
	\put(1.7,0){(a)}
	\put(4.9,0){(b)}
	\put(7.9,0){(c)}
	
	\end{picture}    
	\caption{
	\small{
		Samples of the training and test datasets. (a) The stroke is represented by circular inclusion. This type of inclusion is considered for the training dataset and for the first set of the test dataset. (b) A conductivity sample from the second test dataset where the stroke is modeled as an elliptic inclusion. (c) A conductivity sample from the third test dataset where the stroke is modeled as multiple inclusions containing one disc and one elliptic inclusion.
	}
	}
	\label{fig:stroke_models}
\end{figure}

In real-world scenarios, it is advantageous if a fairly simplistic dataset can be used to train a network which is then made to generalize to more diverse and complicated scenarios. For our training dataset, therefore, the inclusion is always simulated as a disc, see Fig.~\ref{fig:stroke_models}(a). For each sample, the center of the disc and its radius are randomly generated. The radius is randomly
sampled from a uniform distribution in the range $[0.1, 0.25]$, that is ${\cal U}(0.1,0.25)$. Additionally, the conductivity value of the inclusion is randomly sampled from a Gaussian distribution with mean and standard deviation in the ranges given in Table~\ref{table:conductivities}. 
The training dataset contains $4000$ samples, and in order to have a balanced dataset, approximately $50\%$ of the training data has a disc inclusion representing a hemorrhagic stroke and $50\%$ a disc inclusion representing an ischemic stroke. 

To test the performance of the proposed classification method we consider three different and independent test sets where the stroke is modeled as:
(i) one circular inclusion generated in the same way as the training data, see Fig.~\ref{fig:stroke_models}(a);
(ii) one elliptic inclusion with semi-major axis uniformly drawn from the distribution ${\cal U}(0.16,0.26)$, semi-minor axis uniformly drawn from the distribution ${\cal U}(0.07,0.16)$ and randomly rotated, see 
Fig.~\ref{fig:stroke_models}(b); and
(iii) two disjoint inclusions consisting of one disc inclusion and one elliptic inclusion (both representing hemorrhagic strokes or both representing ischemic strokes). In this last case the radius of the circular inclusion is uniformly drawn from the distribution ${\cal U}(0.1,0.15)$, and the semi-major and semi-minor axis of the elliptic inclusion are uniformly drawn from the distributions ${\cal U}(0.16,0.26)$ and  ${\cal U}(0.07,0.16)$, respectively. The elliptic inclusion is also randomly rotated, see Fig.~\ref{fig:stroke_models}(c).

As with the training data, the shape, location and conductivity value of each compartment of the head model is randomly sampled to create the test data. Each test set consists of 1000 samples and  approximately $50\%$ of the samples simulate a hemorrhagic stroke and $50 \%$ an ischemic stroke.

\section{Results and discussion}\label{sec:resuls-discussion}

We present the results of classifying strokes as ischemic or hemorrhagic using two different types of input data: raw EIT voltage measurements and VHED functions. As we mention in Section~\ref{sec:training-testing},  to train the FCNN we use a dataset having only circular inclusions. However, to determine whether the network can generalize to more diverse scenarios, we test the trained network with three different datasets: circular inclusions, elliptical inclusions, and multiple inclusions.

To measure the performance of the FCNN we use the following metrics: (i) \textit{sensitivity --} proportion of actual positives (hemorrhagic cases) that are correctly classified as such, (ii) \textit{specificity --} proportion of actual negatives (ischemic cases) that are correctly identified as such, and (iii) \textit{accuracy --} fraction of correctly classified cases. For each of 20 network trainings, these metrics were computed  on the testing dataset, and the values reported correspond to the average values computed over the 20 trainings. Moreover, during each training process the data was randomly shuffled and the initial weights and biases were chosen randomly.

\subsection{Classification performances}

In Table~\ref{table:circular_inclusion} we present the performance metrics of the FCNN that was trained using a set of 4000 samples where the stroke is modeled as a disc inclusion and tested considering an independent set of 1000 samples also having a disc inclusion. 
We present the results when we assume that data is free of noise and when data is corrupted by two different levels of noise, namely $\delta=10^{-3}$ and $\delta=10^{-2}$. The higher noise level of $\delta =10^{-2}$ may be considered the most realistic scenario. We have highlighted in bold the highest accuracy value for each case.

{\small
	\ctable[
	caption = Performance metrics of FCNN evaluated on a test set of 1000 circular inclusions.,
	label = table:circular_inclusion,
	pos = h
	]{lcccc}{}{
		\FL
		Noise level & Input & Accuracy & Sensitivity & Specificity \ML
		\multirow{2}{*}{$\delta=0$}
		& $V$                       & {\bf 0.987} & 0.988 & 0.987 \NN
		& $T_{\mbox{\tiny odd}}\mu$ & 0.985       & 0.980 & 0.990 \ML
		\multirow{2}{*}{$\delta=10^{-3}$}
		& $V^{\delta}$               &  0.971      & 0.964 & 0.978 \NN
		& $T_{\mbox{\tiny odd}}\mu^{\delta}$ & {\bf 0.984} & 0.976 & 0.992  \ML
		\multirow{2}{*}{$\delta=10^{-2}$}
		& $V^{\delta}$               & 0.937 & 0.935  & 0.937  \NN
		& $T_{\mbox{\tiny odd}}\mu^{\delta}$ & {\bf 0.961} & 0.953  & 0.971    \LL
}}

Table~\ref{table:elliptic_inclusion} shows the performance metrics of the FCNN that was trained using the same set of 4000 samples having a disc inclusion, but was then tested using a set of 1000 samples where the stroke is modeled as an elliptic inclusion, see Fig.~\ref{fig:stroke_models}(b).  
Since the elliptic inclusions differ in shape from the training dataset, this case study serves to evaluate the ability of neural networks to generalize to inputs completely different from those used during training.

{\small
	\ctable[
	caption = Performance metrics of FCNN evaluated on a test set of 1000 elliptic inclusions.,
	label = table:elliptic_inclusion,
	pos = h
	]{lcccc}{}{
		\FL
		Noise level & Input & Accuracy & Sensitivity & Specificity \ML
		\multirow{2}{*}{$\delta=0$}
		& $V$                       & {\bf 0.969} & 0.957 & 0.982 \NN
		& $T_{\mbox{\tiny odd}}\mu$ & 0.965 & 0.943 & 0.987\ML
		\multirow{2}{*}{$\delta=10^{-3}$}
		& $V^{\delta}$               & 0.947 & 0.925 & 0.968 \NN
		& $T_{\mbox{\tiny odd}}\mu^{\delta}$ & {\bf 0.962} & 0.940 & 0.986\ML
		\multirow{2}{*}{$\delta=10^{-2}$}
		& $V^{\delta}$               & 0.885 & 0.865  & 0.905  \NN
		& $T_{\mbox{\tiny odd}}\mu^{\delta}$ & {\bf 0.917} & 0.889  & 0.947  \LL
}}

In Table~\ref{table:multiple_inclusions} we report
the results of the classification problem when the FCNN was trained using a set of 4000 samples having a disc inclusion and evaluated on a set of 1000 samples where the stroke is modeled as two disjoint inclusions shaped like a disc and ellipse, respectively, see Fig.~\ref{fig:stroke_models}(c).
Note that in this case the test dataset differs from the training dataset not only in shape, but also in the number of inclusions.

{\small
	\ctable[
	caption = Performance metrics of FCNN evaluated on a test set of 1000 multiple inclusions.,
	label = table:multiple_inclusions,
	pos = h
	]{lcccc}{}{
		\FL
		Noise level & Input & Accuracy & Sensitivity & Specificity \ML
		\multirow{2}{*}{$\delta=0$}
		& $V$                       & {\bf 0.995} & 0.996 & 0.994 \NN
		& $T_{\mbox{\tiny odd}}\mu$ & 0.989 & 0.984 & 0.993\ML
		\multirow{2}{*}{$\delta=10^{-3}$}
		& $V^{\delta}$               & 0.989 & 0.990 & 0.988 \NN
		& $T_{\mbox{\tiny odd}}\mu^{\delta}$ & 0.989 & 0.982 & 0.995 \ML
		\multirow{2}{*}{$\delta=10^{-2}$}
		& $V^{\delta}$               &  0.951 & 0.954  & 0.947 \NN
		& $T_{\mbox{\tiny odd}}\mu^{\delta}$ &  {\bf 0.971} & 0.962  & 0.979 \LL
}}

\subsection{Discussion on the results}

The numerical results presented in the previous section show that both  raw EIT voltage data and  VHED functions are excellent inputs for the stroke classification problem, even in the presence of the resistive skull layer and the shunting effects induced by the CSF. In the ideal noise-free case, use of raw voltages as inputs to the FCNN produced slightly better accuracy values as compared to the use of VHED functions. This is to be expected, given that the computations used to produce the VHED functions will introduce a small amount of unavoidable numerical roundoff and truncation error, effectively adding some noise to the inputs. However, in the more realistic test scenarios where the simulated voltages were intentionally contaminated with added noise, the VHED functions performed much better than the raw voltages, and the difference is more evident when the noise level is higher (observe the more realistic case $\delta=10^{-2}$). Despite the expected performance decrease for both methods in the presence of increased noise, the accuracy in the case of the VHED functions remains above $91\%$ for all three different models of stroke inclusions, while in the case of raw data inputs, the accuracy drops to $88\%$ for elliptic inclusions. These results reinforce the computational evidence presented in~\cite{agnelli2020classification} which strongly suggests that VHED functions are more robust against noise corruption as compared to raw voltages.

The computation of the VHED functions can be interpreted as a non-linear data preprocessing step for the extraction of noise-robust geometric information. While these VHED functions may contain useful information in their own right, the point of this study was to establish whether VHED functions are useful as inputs for machine learning classification of strokes when applied to physically detailed and mathematically realistic models. Furthermore, our results demonstrate that a simple FCNN trained using single circular-shaped strokes was able to generalize to data generated from phantoms where stroke was modeled as an elliptic inclusion, and even in the case of multiple strokes modeled as a disc and an elliptic inclusion. This provides evidence that relatively simple training data may be used with great success, as networks trained using simplified training data will easily generalize to more diverse scenarios. 
In Table~\ref{table:elliptic_inclusion}, we can see that in the case of elliptic-shaped strokes there is a drop in accuracy for both methods as compared to the case of disc inclusions, which is to be expected given that the training phantoms differ from the test phantoms in this case. In contrast, in Table~\ref{table:multiple_inclusions} we can see that in the case of multiple strokes there is an improvement in accuracy for both methods, which is likely explained by the large physical size of the total stroke region for this scenario, since the larger total stroke size will correspond to a greater perturbation in boundary measurements. These results are consistent with those obtained in~\cite{agnelli2020classification}, where highly simplistic and mathematically unrealistic models were used. The consistency of these results adds credence to the conclusions of this present study, and suggests that the methods presented here are likely to be robust in a variety of different physical scenarios. 

Observe that in all three cases the specificity for the networks using the VHED functions is over $94\%$ for all noise levels. Furthermore, in all scenarios except the noise-free case for multiple strokes, VHED functions yielded higher specificity than raw voltages. Since specificity represents the proportion of actual ischemic cases that are correctly identified, we therefore may conclude that VHED functions are capable of identifying ischemic strokes, especially, with high precision. 
This is important because the majority of all strokes and stroke-related fatalities are ischemic. In fact, as reported in \cite{Grysiewicz2008}, in the United States it has been found that about 87\% of all strokes and 70-80\% of all stroke deaths are ischemic. If the type of stroke can be identified  while still in the ambulance on the way to the hospital, then appropriate medication can be administered immediately, significantly improving the prognosis of the patient.

In this work we employed a very simple FCNN with excellent results. We consider it an advantage of the method that more complicated machine learning architecture did not prove necessary. One remaining simplifying assumption present in this work, however, was the use of a two-dimensional model. 
In practical real-world situations, where a belt or helmet with electrodes is placed on the head of a patient and the EIT data comes from a three-dimensional target, a simple approach like the one used here may not be sufficient. However, for that case we may extend the results presented in this work using more complex architectures or  more sophisticated neural networks, and also employ a more comprehensive set of VHED profiles. This present work therefore serves as a springboard for future work toward clinical application of these methods.

\section{Conclusions}\label{sec:conclusions}

In this article a methodology for classification of stroke as either ischemic or hemorrhagic from \emph{in silico} EIT data was presented and applied to physically detailed and mathematically realistic models which represented a 2D slice of a human head. The phantoms consisted of five regions relevant to electrophysiology of the human head. This included the highly resistive skull layer which impedes current flow, as well as the highly conductive cerebrospinal fluid which introduces shunting effects. Previous studies have frequently omitted one or both of these features from the models used, as the large contrast in the conductivity of these layers makes the classification problem much more difficult. Conductivity values for all tissues and stoke inclusions were furthermore drawn from statistically realistic distributions.    
Following the research line introduced in~\cite{agnelli2020classification}, we extracted noise-robust geometric features known as VHED functions from EIT data, and used these features as inputs for the stroke classification problem. 
In this study we presented a more realistic scenario than in our previous work~\cite{agnelli2020classification}, which employed highly simplified disc-shaped phantoms and used the mathematically ideal continuum model of EIT.  In this work, in addition to the use of a much more physically realistic and detailed phantom, we have considered electrode data simulated using the complete electrode model (CEM) instead of continuous data coming from the idealized continuum model.

We compared the use of VHED functions as network inputs to  the technique employed in other works, which involves learning directly from the raw electrode data. We generated large datasets representing various stroke scenarios that were used to train and test fully connected neural networks using both methods, under three different noise levels. The performance metrics of the classification problem yielded promising  results and demonstrated that in the presence of realistic noise levels, using the VHED functions is a better strategy than directly using EIT voltage measurements in terms of the accuracy and specificity metrics used to evaluate the results. This present work confirms and bolsters the conclusions made in the much more preliminary work~\cite{agnelli2020classification}.

\section*{Acknowledgements}

J.P.A. acknowledges support from ANPCyT project PICT 2021-0188 and SeCyT-UNC Project Formar 2023-33820230100083CB.
S.R. and S.S. acknowledge support from the Research Council of Finland (353097).
M.L. was partially supported by a AdG project 101097198 of the European Research Council, Centre of Excellence of Research Council of Finland and the FAME flagship of the Research Council of Finland (grant 359186). Views and opinions expressed are those of the authors only and do not necessarily reflect those of the funding agencies or the  EU.


\bibliographystyle{plain}
\bibliography{bibliography}
\end{document}